\documentclass[preprint,12pt]{elsarticle}
\usepackage{amsmath}
\usepackage{amssymb}

\journal{}
\input{tcilatex}

\begin{document}

\begin{frontmatter}



\title{New results on equilibria of fuzzy abstract economies}


\author{Monica Patriche}

\address{University of Bucharest, Faculty of Mathematics and Computer Science,  Academiei 14 Street,010014Bucharest,Romania. E-mail:monica.patriche@yahoo.com}
\address{}
\begin{abstract}
We obtain new equilibrium theorems for fuzzy abstract economies with correspondences being w-upper semicontinuous or having e-USS-property. 
\end{abstract}


\end{frontmatter}



\label{}





\bibliographystyle{model1a-num-names}
\bibliography{<your-bib-database>}








\section{Introduction}

Since the theory of fuzzy sets, initiated by Zadeh [23], was considered as a
framework for phenomena which can not be characterized precisely, a lot of
extensions of game theory results have been established. So, many theorems
concerning fuzzy equlibrium existence for fuzzy abstract economies were
obtained. In [10] the authors introduced the concept of a fuzzy game and
proved the existence of equlibrium for 1-person fuzzy game. Also, the
existence of equilibrium points of fuzzy games was studied in [8], [9],
[11],[12], [13], [14], [15], [16], [19]. Fixed point theorems for fuzzy
mappings were proven in [1], [2], [7], [10].

There are several generalizations of the classical model of abstract economy
proposed in his pioneering works by Debreu [4] or later by Shafer and
Sonnenschein [18], Yannelis and Prahbakar [21]. In this paper we consider a
fuzzy extension of Yuan's model of the abstract economy [22] and we prove
the existence of fuzzy equilibrium of fuzzy abstract economies in several
cases. We define two types of correspondences: w-upper semicontinuous
correspondences and correspondences that have e-USS-property. By using a
fixed point theorem for \textit{w-}upper semicontinuous correspondences
[17], we prove our first theorem of equilibrium existence for abstract
economies having w-upper semicontinuous constraint and preference
correspondences. The considered fixed theorem is a Wu like result [20] and
generalizes the Himmelberg's fixed point theorem in [6]. On the other hand,
we use a technique of approximation to prove an equilibrium existence
theorem for set valued maps having e-USS-property.

The paper is organized in the following way: Section 2 contains
preliminaries and notation. The weakly upper semicontinuous correspondences
with respect to a set and the fixed point theorem are presented in Section
3. The equilibrium theorems are stated in Section 4.\bigskip

\section{\textbf{Preliminaries and notation\protect\smallskip \protect%
\medskip }}

Throughout this paper, we shall use the following notations and definitions:

Let $A$ be a subset of a topological space $X$. $\tciFourier (A)$ denotes
the family of all nonempty finite subsets of A. $2^{A}$ denotes the family
of all subsets of $A$. cl$A$ denotes the closure of $A$ in $X$. If $A$ is a
subset of a vector space, co$A$ denotes the convex hull of $A$. If $F$, $G:$ 
$X\rightarrow 2^{Y}$ are correspondences, then co$G$, cl$G$, $G\cap F$ $:$ $%
X\rightarrow 2^{Y}$ are correspondences defined by $($co$G)(x)=$co$G(x)$, $($%
cl$G)(x)=$cl$G(x)$ and $(G\cap F)(x)=G(x)\cap F(x)$ for each $x\in X$,
respectively. The graph of $T:X\rightarrow 2^{Y}$ is the set Gr$%
(T)=\{(x,y)\in X\times Y\mid y\in T(x)\}.$

The correspondence $\overline{T}$ is defined by $\overline{T}(x)=\{y\in
Y:(x,y)\in $cl$_{X\times Y}$Gr$T\}$ (the set cl$_{X\times Y}$Gr$(T)$ is
called the adherence of the graph of $T$)$.$ It is easy to see that cl$%
T(x)\subset \overline{T}(x)$ for each $x\in X.\medskip $

\textbf{Notation.\ }Let E and F be two Hausdorff topological vector spaces
and $X\subset E$, $Y\subset F$ be two nonempty convex subsets. We denote by $%
\mathcal{F}(Y)$ the collection of fuzzy sets on $Y$. A mapping from $X$ into 
$\mathcal{F}(Y)$ is called a fuzzy mapping. If $F:X\rightarrow \mathcal{F}%
(Y) $ is a fuzzy mapping, then for each $x\in X,$ $F(x)$ (denoted by $F_{x}$
in this sequel) is a fuzzy set in $\mathcal{F}(Y)$ and $F_{x}(y)$ is the
degree of membership of point $y$ in $F_{x}.$

A fuzzy mapping $F:X\rightarrow \mathcal{F}(Y)$ is called convex, if for
each $x\in X,$ the fuzzy set $F_{x}$ on $Y$ is a fuzzy convex set, i.e., for
any $y_{1},y_{2}\in Y,$ $t\in \lbrack 0,1],$ $F_{x}(ty_{1}+(1-t)y_{2})\geq
\min \{F_{x}(y_{1}),F_{x}(y_{2})\}.$

In the sequel, we denote by

$(A)_{q}=\{y\in Y:A(y)\geq q\},$ $q\in \lbrack 0,1]$ the $q$-cut set of $%
A\in \mathcal{F}(Y).\medskip $

\begin{definition}
Let $X$, $Y$ be topological spaces and $T:X\rightarrow 2^{Y}$ be a
correspondence. $T$ is said to be \textit{upper semicontinuous} if for each $%
x\in X$ and each open set $V$ in $Y$ with $T(x)\subset V$, there exists an
open neighborhood $U$ of $x$ in $X$ such that $T(y)\subset V$ for each $y\in
U$. $T$ is said to be \textit{almost upper semicontinuous} if for each $x\in
X$ and each open set $V$ in $Y$ with $T(x)\subset V$, there exists an open
neighborhood $U$ of $x$ in $X$ such that $T(y)\subset $cl$V$ for each $y\in
U $.\medskip
\end{definition}

\begin{lemma}
(Lemma 3.2, pag. 94 in [24])Let $X$ be a topological space, $Y$ be a
topological linear space, and let $S:X\rightarrow 2^{Y}$ be an upper
semicontinuous correspondence with compact values. Assume that the sets $%
C\subset Y$ and $K\subset Y$ are closed and respectively compact. Then $%
T:X\rightarrow 2^{Y}$ defined by $T(x)=(S(x)+C)\cap K$ for all $x\in X$ is
upper semicontinuous.\medskip
\end{lemma}

Lemma 2 is a version of Lemma 1.1 in [22] ( for $D=Y,$ we obtain Lemma 1.1
in [22]). For the reader's convenience, we include its proof below.

\begin{lemma}
Let $X$ be a topological space, $Y$ be a nonempty subset of a locally convex
topological vector space $E$ and $T:X\rightarrow 2^{Y}$ be a correspondence$%
. $ Let \ss\ be a basis of neighbourhoods of $0$ in $E$ consisting of open
absolutely convex symmetric sets. Let $D$ be a compact subset of $Y$. If for
each $V\in $\ss , the correspondence $T^{V}:X\rightarrow 2^{Y}$ is defined
by $T^{V}(x)=(T(x)+V)\cap D$ for each $x\in X,$ then $\cap _{V\in \text{\ss }%
}\overline{T^{V}}(x)\subseteq \overline{T}(x)$ for every $x\in X.\medskip $
\end{lemma}

\textit{Proof.} Let $x$ and $y$ be such that $y\in \cap _{V\in \text{\ss }}%
\overline{T^{V}}(x)$ and suppose, by way of contradiction, that $y\notin 
\overline{T}(x).$ This means that $(x,y)\notin $cl Gr $T,$ so that there
exists an open neighborhood $U$ of $x$ and $V\in $\ss\ such that:

$(U\times (y+V))\cap $Gr $T=\emptyset .\ \ \ \ \ \ \ \ \ \ \ \ \ \ \ \ \ \ \
\ \ \ \ \ \ \ \ \ \ \ \ \ (1)$

Choose $W\in $\ss\ such that $W-W\subseteq V$ (e.g. $W=\frac{1}{2}V)$. Since 
$y\in T^{W}(x)$, then $(x,y)\in $cl Gr $T^{W},$ so that

\begin{equation*}
(U\times (y+W))\cap \text{Gr }T^{W}\neq \emptyset .\ \ \ \ \ \ \ \ \ \ \ \ \
\ \ \ \ \ \ \ \ \ \ \ \ \ \ \ \ \ \ \ 
\end{equation*}

There are some $x^{\prime }\in U$ and $w^{\prime }\in W$ such that $%
(x^{\prime },y+w^{\prime })\in $Gr $T^{W},$ i.e. $y+w^{\prime }\in
T^{W}(x^{\prime }).$ Then, $y+w^{\prime }\in D$ and $y+w^{\prime }=y^{\prime
}+w^{^{\prime \prime }}$ for some $y^{\prime }\in T(x^{\prime })$ and $%
w^{^{\prime \prime }}\in W.$ Hence, $y^{\prime }=y+(w^{\prime }-w^{^{\prime
\prime }})\in y+(W-W)\subseteq y+V,$ so that $T(x^{\prime })\cap (y+V)\neq
\emptyset .$ Since $x^{\prime }\in U,$ this means that $(U\times (y+V))\cap $%
Gr $T\neq \emptyset ,$ contradicting (1). $\ \ \ \ \ \ \ \ \ \ \square $%
\medskip

\section{Weakly upper semicontinuous correspondences with respect to a set}

We introduce the following definitions.

Let $X$ be a topological space, $Y$ be a nonempty subset of a topological
vector space $E$ and $D$ be a subset of $Y$.

\begin{definition}
The correspondence $T:X\rightarrow 2^{Y}$ is said to be \textit{w-upper
semicontinuous} (weakly upper semicontinuous) \textit{with respect to} 
\textit{the set} $D$ if there exists a basis \ss\ of open symmetric
neighborhoods of $0$ in $E$ such that, for each $V\in $\ss , the
correspondence $T^{V}$ is upper semicontinuous.
\end{definition}

\begin{definition}
The correspondence $T:X\rightarrow 2^{Y}$ is said to be \textit{almost} 
\textit{w-upper semicontinuous} (almost weakly upper semicontinuous) \textit{%
with respect to} \textit{the set} $D$ if there exists a basis \ss\ of open
symmetric neighborhoods of $0$ in $E$ such that, for each $V\in $\ss , the
correspondence $\overline{T^{V}}$ is upper semicontinuous.\medskip
\end{definition}

\begin{example}
Let $T_{1}:(0,2)\rightarrow 2^{(0,2)}$ be defined by $T_{1}(x)=\left\{ 
\begin{array}{c}
(0,1)\text{ if }x\in (0,1]; \\ 
\lbrack 1,2)\text{ if x}\in (1,2).%
\end{array}%
\right. $
\end{example}

$T_{1}$ and $T_{1}\cap \{1\}=\left\{ 
\begin{array}{c}
\phi \text{ \ if \ }x\in (0,1]; \\ 
\{1\}\text{ if x}\in (1,2)%
\end{array}%
\right. $ are not upper semicontinuous on $(0,2),$ but $T_{1}$ is w-upper
semicontinuous with respect to $D$ and it is also almost w-upper
semicontinuous with respect to $D.$

We also define the dual w-upper semicontinuity with respect to a compact set.

\begin{definition}
Let $T_{1},T_{2}:X\rightarrow 2^{Y}$ be correspondences. The pair $%
(T_{1},T_{2})$ is said to be \textit{dual almost w-upper semicontinuous}
(dual almost weakly upper semicontinuous) \textit{with respect to} the set $%
D $ if there exists a basis \ss\ of open symmetric neighborhoods of $0$ in $%
E $ such that, for each $V\in $\ss , the correspondence $\overline{%
T_{(1,2)}^{V}}:X\rightarrow 2^{D}$ is lower semicontinuous, where $%
T_{(1,2)}^{V}:X\rightarrow 2^{D}$ is defined by $%
T_{(1,2)}^{V}(x)=(T_{1}(x)+V)\cap T_{2}(x)\cap D$ for each $x\in X$.\medskip
\end{definition}

\begin{example}
Let $\ D=[1,2],$ $T_{1}:(0,2)\rightarrow 2^{[1,4]}$ be the correspondence
defined by
\end{example}

$T_{1}(x)=\left\{ 
\begin{array}{c}
\lbrack 2-x,2],\text{ if }x\in (0,1); \\ 
\{4\}\text{ \ \ \ \ \ \ if \ \ \ \ \ \ \ }x=1; \\ 
\lbrack 1,2]\text{ \ \ \ if \ \ \ }x\in (1,2).%
\end{array}%
\right. $

and $T_{2}:(0,2)\rightarrow 2^{[2,3]}$ be the correspondence defined by

$T_{2}(x)=\left\{ 
\begin{array}{c}
\lbrack 2,3],\text{ if }x\in (0,1]; \\ 
\{2\}\text{ \ if \ \ }x\in (1,2);%
\end{array}%
\right. .$

The correspondence $T_{1}$ is not upper semicontinuous on $(0,2)$, but $%
\overline{T_{(1,2)}^{V}}$ is upper semicontinuous and has nonempty values.

We conclude that the pair $(T_{1},T_{2})$ is dual almost w-upper
semicontinuous with respect to $D.$

We obtain the following fixed point theorem which generalizes Himmelberg's
fixed point theorem in [6]:

\begin{theorem}
(see [17])\textit{Let }$I$\textit{\ be an index set. For each }$i\in I,$%
\textit{\ let }$X_{i}$\textit{\ be a nonempty convex subset of a Hausdorff
locally convex topological vector space }$E_{i}$\textit{, }$D_{i}$\ be a
nonempty compact convex subset of $X_{i}$ \textit{and }$S_{i},T_{i}:X:=%
\tprod\limits_{i\in I}X_{i}\rightarrow 2^{X_{i}}$\textit{\ be two
correspondences with the following conditions:}
\end{theorem}

1) \textit{for each }$x\in X$, $\overline{S}_{i}(x)\subseteq T_{i}(x)$%
\textit{. }

2)\textit{\ }$S_{i}$\textit{\ is almost w-upper semicontinuous with respect
to }$D_{i}$ \textit{and} $\overline{S_{i}^{V_{i}}}$ \textit{is} \textit{%
convex nonempty valued for each absolutely convex symmetric neighborhood }$%
V_{i}$ \textit{of} $0$ \textit{in} $E_{i}$\textit{.}

\textit{Then there exists} $x^{\ast }\in D:=\tprod\limits_{i\in I}D_{i}$%
\textit{\ such that }$x_{i}^{\ast }\in T_{i}(x^{\ast })$ \textit{for each} $%
i\in I.\medskip $

\section{Existence of fuzzy equilibrium for fuzzy abstract economies\textit{%
\ }}

\subsection{The model of a fuzzy abstract economy}

In this section we describe the fuzzy equilibrium for a fuzzy extension of
Yuan's model of abstract economy [22]. We prove the existence of fuzzy
equilibrium of abstract fuzzy economies in several cases.

Let $I$ be a nonempty set (the set of agents). For each $i\in I$, let $X_{i}$
be a non-empty topological vector space representing the set of actions and
define $X:=\underset{i\in I}{\prod }X_{i}$; let $A_{i}$, $B_{i}:X\rightarrow 
$ $\mathcal{F}(X_{i})$ be the constraint fuzzy correspondences and $P_{i}$ $%
:X\rightarrow $ $\mathcal{F}(X_{i})$ the preference fuzzy correspondence, $%
a_{i},b_{i}:X\rightarrow (0,1]$ fuzzy constraint functions and $%
p_{i}:X\rightarrow (0,1]$ fuzzy preference function.

Let denote $A_{i}^{\prime }$, $B_{i}^{\prime },P_{i}^{\prime }:X\rightarrow $
$2^{X_{i}},$ defined by $A_{i}^{\prime }(x)=(A_{i_{x}})_{a_{i}(x)}$, $%
B_{i}^{\prime }(x)=(B_{i_{x}})_{b_{i}(x)}$ and $P_{i}^{\prime
}(x)=(P_{i_{x}})_{b_{i}(x)}.$

\textit{Definition} 4.\textit{\ }A fuzzy \textit{abstract economy} is
defined as an ordered family $\Gamma
=(X_{i},A_{i},B_{i},P_{i},a_{i},b_{i},p_{i})_{i\in I}$.

If $A_{i},B_{i},P_{i}:X\rightarrow 2^{Y_{i}}$ are classical correspondences,
then the previous definition can be reduced to the standard definition of
abstract\textit{\ }economy due to Yuan [22].\medskip

\textit{Definition} 5.\textit{\ }A \textit{fuzzy} \textit{equilibrium} for $%
\Gamma $ is defined as a of point $x^{\ast }\in X$ such that for each $i\in
I $, $x_{i}^{\ast }\in \overline{B^{\prime }}(x^{\ast })$ and $%
(A_{i_{x^{\ast }}})_{a_{i}(x^{\ast })}\cap (P_{i_{x^{\ast
}}})_{p_{i}(x^{\ast })}=\emptyset ,$ where $(A_{i_{x^{\ast
}}})_{a_{i}(x^{\ast })}=\{z\in Y_{i}:A_{i_{x^{\ast }}}(z)\geq a_{i}(x^{\ast
})\},$ $(B_{i_{x^{\ast }}})_{b_{i}(x^{\ast })}=\{z\in Y_{i}:B_{i_{x^{\ast
}}}(z)\geq b_{i}(x^{\ast })\}$, $(P_{i_{x^{\ast }}})_{p_{i}(x^{\ast
})}=\{z\in Y_{i}:P_{i_{x^{\ast }}}(z)\geq p_{i}(x^{\ast })\}$.\medskip

\subsection{Equilibria existence}

As an application of the fixed point Theorem 1, we have the following
result.\medskip

\begin{theorem}
\textit{Let }$\Gamma =(X_{i},A_{i},B_{i},P_{i},a_{i},b_{i},p_{i})_{i\in I}$%
\textit{\ be a fuzzy abstract economy such that for each }$i\in I$\textit{\
the following conditions are fulfilled:}
\end{theorem}

\textit{1)\ }$X_{i}$\textit{\ \ be a non-empty compact convex subset of a
locally convex Hausdorff topological vector space }$E_{i}$\textit{\ and }$%
D_{i}$\textit{\ is a nonempty compact convex subset of }$X_{i}$\textit{;}

\textit{2)\ }$A_{i},P_{i}$\textit{\ and }$B_{i}$\textit{\ are such that each 
}$(B_{i_{x}})_{b_{i}(x)}$\textit{\ is a nonempty convex subset of }$X_{i},$%
\textit{\ }$(A_{i_{x}})_{a_{i}(x)},(P_{i_{x}})_{p_{i}(x)}$\textit{\ are
convex and }$(A_{i_{x}})_{a_{i}(x)}\cap (P_{i_{x}})_{p_{i}(x)}\subset
(B_{i_{x}})_{b_{i}(x)}$\textit{\ for each }$x\in X;$

\textit{3) the set }$W_{i}=\left\{ x\in X:(A_{i_{x}})_{a_{i}(x)}\cap
(P_{i_{x}})_{p_{i}(x)}\neq \emptyset \right\} $\textit{\ is open in }$X$%
\textit{.}

\textit{4) the correspondence }$H_{i}:X\rightarrow 2^{X_{i}}$\textit{\
defined by }$H_{i}\left( x\right) =(A_{i_{x}})_{a_{i}(x)}\cap
(P_{i_{x}})_{p_{i}(x)}$\textit{\ for each }$x\in X$\textit{\ is almost
w-upper semicontinuous with respect to }$D_{i}$\textit{\ on }$W_{i}$\textit{%
\ and }$\overline{H_{i}^{V_{i}}}$\textit{\ is convex nonempty valued for
each open absolutely convex symmetric neighborhood }$V_{i}$\textit{\ of }$0$%
\textit{\ in }$E_{i}$\textit{;}

\textit{5)\ the correspondence }$x\rightarrow
(B_{i_{x}})_{b_{i}(x)}:X\rightarrow 2^{X_{i}}$\textit{\ is almost w-upper
semicontinuous with respect to }$D_{i}$\textit{\ and }$\overline{%
B_{i}^{V_{i}}}$\textit{\ is convex nonempty valued for each open absolutely
convex symmetric neighborhood }$V_{i}$\textit{\ of }$0$\textit{\ in }$E_{i},$
where $B_{i}^{V_{i}}:X\rightarrow 2^{X_{i}}$ is defined by $%
B_{i}^{V_{i}}(x)=((B_{i_{x}})_{b_{i}(x)}+V_{i})\cap D_{i}$\textit{;}

\textit{6)\ for each }$x\in X$\textit{\ , }$x_{i}\notin \overline{H_{i}}(x)$%
\textit{;}

\textit{Then there exists a fuzzy equilibrium point }$x^{\ast }\in D=\
\prod\limits_{i\in I}D_{i}$\textit{\ such that for each }$i\in I$\textit{, }$%
x_{i}^{\ast }\in \overline{B}_{i}^{\prime }(x^{\ast })$\textit{\ and }$%
(A_{i_{x^{\ast }}})_{a_{i}(x^{\ast })}\cap (P_{i_{x^{\ast
}}})_{p_{i}(x^{\ast })}=\emptyset .$\textit{\medskip }

\textit{Proof. }Let $i\in I.$ By condition (3) we know that $W_{i}$ is open
in $X.$

Let's define $T_{i}:X\rightarrow 2^{X_{i}}$ by $T_{i}\left( x\right)
=\left\{ 
\begin{array}{c}
(A_{i_{x}})_{a_{i}(x)}\cap (P_{i_{x}})_{p_{i}(x)},\text{ if }x\in W_{i}, \\ 
(B_{i_{x}})_{b_{i}(x)},\text{ \ \ \ \ \ \ \ \ \ \ \ if }x\notin W_{i}%
\end{array}%
\right. $ for each $x\in X.$

Then $T_{i}:X\rightarrow 2^{X_{i}}$ is a correspondence with nonempty convex
values. We shall prove that $T_{i}:X\rightarrow 2^{D_{i}}$ is almost w-upper
semicontinuous with respect to $D_{i}$. Let \ss $_{i}$\ be a basis of open
absolutely convex symmetric neighborhoods of $0$ in $E_{i}$ and let \ss =$%
\tprod\limits_{i\in I}$\ss $_{i}.$

For each $V=(V_{i})_{i\in I}\in \tprod\limits_{i\in I}$\ss $_{i},$ for each $%
x\in X,$ let for each $i\in I$

$F_{i}^{V_{i}}(x)=((A_{i_{x}})_{a_{i}(x)}\cap
(P_{i_{x}})_{p_{i}(x)}+V_{i})\cap D_{i}$ and

$T_{i}^{V_{i}}(x)=\left\{ 
\begin{array}{c}
F_{i}^{V_{i}}(x),\text{ if }x\in W_{i}, \\ 
B_{i}^{V_{i}}(x),\text{\ if }x\notin W_{i}.%
\end{array}%
\right. $

For each open set $V_{i}^{\prime }$ in $D_{i}$, the set

$\left\{ x\in X:\overline{T_{i}^{V_{i}}}\left( x\right) \subset
V_{i}^{\prime }\right\} =$

$=\left\{ x\in W_{i}:\overline{F_{i}^{V_{i}}}(x)\subset V_{i}^{\prime
}\right\} \cup \left\{ x\in X\smallsetminus W_{i}:\overline{B_{i}^{V_{i}}}%
(x)\subset V_{i}^{^{\prime }}\right\} $

$=\left\{ x\in W_{i}:\overline{F_{i}^{V_{i}}}(x)\subset V_{i}^{^{\prime
}}\right\} \cup \left\{ x\in X:\overline{B_{i}^{V_{i}}}(x)\subset
V_{i}^{\prime }\right\} .$

According to condition (4), the set $\left\{ x\in W_{i}:\overline{%
F_{i}^{V_{i}}}(x)\subset V_{i}^{\prime }\right\} $ is open in $X$. The set $%
\left\{ x\in X:\overline{B_{i}^{V_{i}}}(x)\subset V_{i}^{\prime }\right\} $
is open in $X$ because $\overline{B_{i}^{V_{i}}}$ is upper semicontinuous.

Therefore, the set $\left\{ x\in X:\overline{T_{i}^{V_{i}}}\left( x\right)
\subset V_{i}^{\prime }\right\} $ is open in $X.$ It shows that $\overline{%
T_{i}^{V_{i}}}:X\rightarrow 2^{D_{i}}$ is upper semicontinuous. According to
Theorem 1, there exists $x^{\ast }\in D=$ $\prod\limits_{i\in I}D_{i}$ such
that $x^{\ast }\in \overline{T}_{i}\left( x^{\ast }\right) ,$ for each $i\in
I.$ By condition (5) we have that $x_{i}^{\ast }\in \overline{B_{i}^{\prime }%
}(x^{\ast })$ and $(A_{i_{x^{\ast }}})_{a_{i}(x^{\ast })}\cap (P_{i_{x^{\ast
}}})_{p_{i}(x^{\ast })}=\emptyset $ for each $i\in I.\square \medskip $

Theorem 5 deals with abstract economies which have dual w-upper
semicontinuous pairs of correspondences.

\begin{theorem}
\textit{Let }$\Gamma =(X_{i},A_{i},B_{i},P_{i},a_{i},b_{i},p_{i})_{i\in I}$%
\textit{\ be a fuzzy abstract economy such that for each }$i\in I$\textit{\
the following conditions are fulfilled:}
\end{theorem}

\textit{1)\ }$X_{i}$\textit{\ \ be a non-empty compact convex subset of a
locally convex Hausdorff topological vector space }$E_{i}$\textit{\ and }$%
D_{i}$\textit{\ is a nonempty compact convex subset of }$X_{i}$\textit{;}

\textit{2)\ }$A_{i},P_{i},B_{i}$\textit{\ are such that each }$%
(B_{i_{x}})_{b_{i}(x)}$\textit{\ is a convex subset of }$X_{i},$\textit{\ }$%
(P_{i_{x}})_{p_{i}(x)}\subset D_{i}$\textit{\ and }$(A_{i_{x}})_{a_{i}(x)}%
\cap (P_{i_{x}})_{p_{i}(x)}\subset (B_{i_{x}})_{b_{i}(x)}$\textit{\ for each 
}$x\in X;$

\textit{3) the set }$W_{i}=\left\{ x\in X:(A_{i_{x}})_{a_{i}(x)}\cap
(P_{i_{x}})_{p_{i}(x)}\neq \emptyset \right\} $\textit{\ is open in }$X$%
\textit{.}

\textit{4) the pair }$x\rightarrow ((A_{i_{x}})_{_{a_{i}(x)}\mid \text{cl}%
W_{i}},(P_{i_{x}})_{_{p_{i}(x)}\mid \text{cl}W_{i}})$\textit{\ is dual
almost w-upper semicontinuous with respect to }$D_{i}$\textit{, the
correspondence }$x\rightarrow (B_{i_{x}})_{_{b_{i}(x)}}:X\rightarrow
2^{X_{i}}$\textit{\ is almost w-upper semicontinuous with respect to }$D_{i}$%
\textit{;}

\textit{5) if }$T_{i,V_{i}}:X\rightarrow 2^{X_{i}}$\textit{\ is defined by }$%
T_{i,V_{i}}(x)=((A_{i_{x}})_{a_{i}(x)}+V_{i})\cap D_{i}\cap
(P_{i_{x}})_{p_{i}(x)}$\textit{\ for each }$x\in X$\textit{\ and }$%
B_{i}^{V_{i}}:X\rightarrow 2^{X_{i}}$\textit{\ is defined by }$%
B_{i}^{V_{i}}(x)=((B_{i_{x}})_{b_{i}(x)}+V_{i})\cap D_{i}$\textit{\ for each 
}$x\in X,$\textit{\ then the correspondences }$\overline{B_{i}^{V_{i}}}$%
\textit{\ and }$\overline{T_{i,V_{i}}}$\textit{\ are nonempty convex valued
for each open absolutely convex symmetric neighborhood }$V_{i}$\textit{\ of }%
$0$\textit{\ in }$E_{i}$\textit{;}

\textit{6)\ for each }$x\in X$\textit{, }$x_{i}\notin \overline{%
P_{i}^{\prime }}(x).$

\textit{Then there exists a fuzzy equilibrium point }$x^{\ast }\in D=\
\prod\limits_{i\in I}D_{i}$\textit{\ such that for each }$i\in I$\textit{, }$%
x_{i}^{\ast }\in \overline{B_{i}^{\prime }}(x^{\ast })$\textit{\ and }$%
(A_{i_{x^{\ast }}})_{a_{i}(x^{\ast })}\cap (P_{i_{x^{\ast
}}})_{p_{i}(x^{\ast })}=\emptyset .$\textit{\medskip }

\textit{Proof.} For each $i\in I,$ let \ss $_{i}$\ denote the family of all
open absolutely convex symmetric neighborhoods of zero in $E_{i}$ and let 
\ss $=\tprod\limits_{i\in I}$\ss $_{i}.$ For each $V=\tprod\limits_{i\in
I}V_{i}\in \tprod\limits_{i\in I}$\ss $_{i},$ for each $i\in I,$ let

$A_{i}^{V_{i}},S_{i}^{V_{i}}:X\rightarrow 2^{X_{i}}$ be defined by

$A_{i}^{V_{i}}(x)=((A_{i_{x}})_{_{b_{i}(x)}}+V_{i})\cap D_{i}$ for each $%
x\in X$ and

$S_{i}^{V_{i}}\left( x\right) =\left\{ 
\begin{array}{c}
T_{i,V_{i}}(x),\text{ \ \ \ \ \ \ \ \ \ \ \ \ \ \ \ \ \ \ if }x\in W_{i}, \\ 
B_{i}^{V_{i}}(x),\text{ \ \ \ \ \ \ \ \ \ \ \ \ \ \ \ \ if }x\notin W_{i},%
\end{array}%
\right. $

$\overline{S_{i}^{V_{i}}}$ has closed values. Next, we shall prove that $%
\overline{S_{i}^{V_{i}}}:X\rightarrow 2^{D_{i}}$ is upper semicontinuous.

For each open set $V^{\prime }$ in $D_{i}$, the set

$\left\{ x\in X:\overline{S_{i}^{V_{i}}}\left( x\right) \subset V^{\prime
}\right\} =$

$=\left\{ x\in W_{i}:\overline{T_{i,V_{i}}}(x)\subset V^{\prime }\right\}
\cup \left\{ x\in X\smallsetminus W_{i}:\overline{B_{i}^{V_{i}}}(x)\subset
V^{\prime }\right\} $

=$\left\{ x\in W_{i}:\overline{T_{i,V_{i}}}(x)\subset V^{\prime }\right\}
\cup \left\{ x\in X:\overline{B_{i}^{V_{i}}}(x)\subset V^{\prime }\right\} .$

We know that the correspondence $\overline{T_{i,V_{i}}}(x)_{\mid W_{i}}:$ $%
W_{i}\rightarrow 2^{D_{i}}$ is upper semicontinuous. The set $\left\{ x\in
W_{i}:\overline{T_{i,V_{i}}}(x)\subset V^{\prime }\right\} $ is \ open in $X.
$ Since $\overline{B_{i}^{V_{i}}}(x):X\rightarrow 2^{D_{i}}$ is upper
semicontinuous, the set $\{x\in X:\overline{B_{i}^{V_{i}}}(x)\}\subset
V^{\prime }$ is open in $X$ and therefore, the set $\left\{ x\in X:\overline{%
S_{i}^{V_{i}}}\left( x\right) \subset V^{\prime }\right\} $ is open in $X$.
It proves that $\overline{S_{i}^{V_{i}}}:X\rightarrow 2^{D_{i}}$ is upper
semicontinuous. According to Himmelberg's Theorem, applied for the
correspondences $\overline{S_{i}^{V_{i}}},$ there exists a point $%
x_{V}^{\ast }\in D=$ $\prod\limits_{i\in I}D_{i}$ such that $(x_{V}^{\ast
})_{i}\in S_{i}^{V_{i}}\left( x_{V}^{\ast }\right) $ for each $i\in I.$ By
condition (5), we have that $(x_{V}^{\ast })_{i}\notin \overline{%
P_{i}^{\prime }}(x_{V}^{\ast }),$ hence, $(x_{V}^{\ast })_{i}\notin 
\overline{A_{i}^{V_{i}}}\left( x_{V}^{\ast }\right) \cap \overline{%
P_{i}^{\prime }}(x_{V}^{\ast })$. \newline
We also have that clGr$(T_{i,V_{i}})\subseteq $ clGr$(A_{i}^{V_{i}})\cap $%
clGr$P_{i}^{\prime }.$ Then $\overline{T_{i,V_{i}}}(x)\subseteq $ $\overline{%
A_{i}^{V_{i}}}(x)\cap \overline{P_{i}^{\prime }}\left( x\right) $ for each $%
x\in X.$ It follows that $(x_{V}^{\ast })_{i}\notin \overline{T_{i,V_{i}}}%
(x_{V}^{\ast }).$ Therefore, $(x_{V}^{\ast })_{i}\in \overline{B_{i}^{V_{i}}}%
\left( x_{V}^{\ast }\right) .$

For each $V=(V_{i})_{i\in I}\in \tprod\limits_{i\in I}$\ss $_{i},$ let's
define $Q_{V}=\cap _{i\in I}\{x\in D:x\in \overline{B_{i}^{V_{i}}}\left(
x\right) $ and $(A_{i_{x}})_{a_{i}(x)}\cap (P_{i_{x}})_{p_{i}(x)}=\emptyset
\}.$

$Q_{V}$ is nonempty since $x_{V}^{\ast }\in Q_{V},$ and it is a closed
subset of $D$ according to (3). Then, $Q_{V}$ is nonempty and compact.

Let \ss =$\tprod\limits_{i\in I}$\ss $_{i}.$ We prove that the family $%
\{Q_{V}:V\in \text{\ss }\}$ has the finite intersection property.

Let $\{V^{(1)},V^{(2)},...,V^{(n)}\}$ be any finite set of $\text{\ss\ }$and
let $V^{(k)}=\underset{i\in I}{\tprod }V_{i}^{(k)}{}_{i\in I}$, $k=1,...,n.$
For each $i\in I$, let $V_{i}=\underset{k=1}{\overset{n}{\cap }}V_{i}^{(k)}$%
, then $V_{i}\in \text{\ss }_{i};$ thus $V\in \underset{i\in I}{\tprod }%
\text{\ss }_{i}.$ Clearly $Q_{V}\subset \underset{k=1}{\overset{n}{\cap }}%
Q_{V^{(k)}}$ so that $\underset{k=1}{\overset{n}{\cap }}Q_{V^{(k)}}\neq
\emptyset .$

Since $D$ is compact and the family $\{Q_{V}:V\in \text{\ss }\}$ has the
finite intersection property, we have that $\cap \{Q_{V}:V\in \text{\ss }%
\}\neq \emptyset .$ Take any $x^{\ast }\in \cap \{Q_{V}:V\in $\ss $\},$ then
for each $V\in \text{\ss },$

$x^{\ast }\in \cap _{i\in I}\left\{ x^{\ast }\in D:x_{i}^{\ast }\in 
\overline{B_{i}^{V_{i}}}\left( x^{\ast }\right) \text{ and }%
(A_{i_{x}})_{a_{i}(x)}\cap (P_{i_{x}})_{p_{i}(x)}=\emptyset )\right\} .$

Hence, $x_{i}^{\ast }\in \overline{B_{i}^{V_{i}}}\left( x^{\ast }\right) $
for each $V\in $\ss\ and for each $i\in I.$ According to Lemma\emph{\ }2,%
\emph{\ }we have that\emph{\ } $x_{i}^{\ast }\in \overline{(B_{i}^{\prime }}%
)(x^{\ast })$ and $(A_{i_{x^{\ast }}})_{a_{i}(x^{\ast })}\cap (P_{i_{x^{\ast
}}})_{p_{i}(x^{\ast })}=\emptyset $ for each $i\in I.$ $\square $\medskip

We now introduce the following concept, which also generalizes the concept
of lower semicontinuous correspondences.

\begin{definition}
Let $X$ be a non-empty convex subset of a topological linear space $E$, $Y$
be a non-empty set in a topological space and $K\subseteq X\times Y.$
\end{definition}

The correspondence $T:X\times Y\rightarrow 2^{X}$ has the e-USCS-property
(e-upper semicontinuous selection property) on $K,$ if for each absolutely
convex neighborhood $V$ of zero in $E,$ there exists an upper semicontinuous
correspondence with convex values $S^{V}:X\times Y\rightarrow 2^{X}$ such
that $S^{V}(x,y)\subset T(x,y)+V$ and $x\notin $cl$S^{V}(x,y)$ for every $%
(x,y)\in K$.$\medskip $

The following theorem is an equilibrium existence result for economies with
constraint correspondences having e-USCS-property.

\begin{theorem}
\textit{Let }$\Gamma =(X_{i},A_{i},B_{i},P_{i},a_{i},b_{i},p_{i})_{i\in I}$%
\textit{\ be a fuzzy abstract economy such that for each }$i\in I$\textit{\
the following conditions are fulfilled:}
\end{theorem}

\textit{1)\ }$X_{i}$\textit{\ \ be a non-empty compact convex subset of a
locally convex Hausdorff space }$E_{i}$\textit{;}

\textit{(2)\ the correspondence }$x\rightarrow $\textit{cl}$%
(B_{i_{x}})_{b_{i}(x)}:X\rightarrow 2^{X_{i}}$\textit{\ is upper
semicontinuous with non-empty convex values;}

\textit{(3)\ the set }$W_{i}:$\textit{\ }$=\left\{ x\in
X:(A_{i_{x}})_{a_{i}(x)}\cap (P_{i_{x}})_{p_{i}(x)}\neq \emptyset \right\} $%
\textit{\ is open;}

\textit{(3)\ the correspondence }$x\rightarrow $\textit{cl}$%
((A_{i_{x}})_{a_{i}(x)}\cap (P_{i_{x}})_{p_{i}(x)}):X\rightarrow 2^{X_{i}}$%
\textit{\ has \ the e-USCS-property on }$W_{i}$\textit{.}

\textit{Then there exists an equilibrium point }$x^{\ast }\in X$\textit{\ \
for }$\Gamma $\textit{,}$\ i.e.$\textit{, for each }$i\in I$\textit{, }$%
x_{i}^{\ast }\in \overline{B^{\prime }}(x^{\ast })$\textit{\ and }$%
(A_{i_{x^{\ast }}})_{a_{i}(x^{\ast })}\cap (P_{i_{x^{\ast
}}})_{p_{i}(x^{\ast })}=\emptyset .$\textit{\medskip }

\textit{Proof.} For each\textit{\ }$i\in I$, let \ss $_{i}$ denote the
family of all open convex neighborhoods of zero in $E_{i}.$ Let $%
V=(V_{i})_{i\in I}\in \tprod\limits_{i\in I}$\ss $_{i}.$ Since the
correspondence $x\rightarrow $cl$((A_{i_{x}})_{a_{i}(x)}\cap
(P_{i_{x}})_{p_{i}(x)})$ has the e-USCS-property on $W_{i}$, it follows that
there exists an upper semicontinuous correspondence $F_{i}^{V_{i}}:X%
\rightarrow 2^{X_{i}}$ such that $F_{i}^{V_{i}}(x)\subset $cl$%
((A_{i_{x}})_{a_{i}(x)}\cap (P_{i_{x}})_{p_{i}(x)})+V_{i}$ and $x_{i}\notin $%
cl$F_{i}^{V_{i}}(x)$ for each $x\in W_{i}$.

Define the correspondence $T_{i}^{V_{i}}:X\rightarrow 2^{X_{i}}$, by

$T_{i}^{V_{i}}(x):=\left\{ 
\begin{array}{c}
\text{cl}\{F_{i}^{V_{i}}(x)\}\text{, \ \ \ \ \ \ \ \ \ \ \ \ \ \ \ \ \ \ if }%
x\in W_{i}\text{, } \\ 
\text{cl}((B_{i_{x}})_{b_{i}(x)}+V_{i})\cap X_{i}\text{, if }x\notin W_{i}%
\text{;}%
\end{array}%
\right. $

$B_{i}^{V_{i}}:X\rightarrow 2^{X_{i}},$ $B_{i}^{V_{i}}(x)=$cl$%
((B_{i_{x}})_{b_{i}(x)}+V_{i})\cap X_{i}=($cl$(B_{i_{x}})_{b_{i}(x)}+$cl$%
V_{i})\cap X_{i}$ is upper semicontinuous by Lemma 1\emph{.}

Let $U$ be an open subset of $\ X_{i}$, then

$U^{^{\prime }}:=\{x\in X$ $\mid T_{i}^{V_{i}}(x)\subset U\}$

\ \ \ =$\{x\in W_{i}$ $\mid T_{i}^{V_{i}}(x)\subset U\}\cup \{x\in
X\setminus W_{i}$ $\mid $ $T_{i}^{V_{i}}(x)\subset U\}$

\ \ \ =$\left\{ x\in W_{i}\text{ }\mid \text{cl}F_{i}^{V_{i}}(x)\subset
U\right\} \cup \left\{ x\in X\mid \text{ }(\text{cl}(B_{i_{x}})_{b_{i}(x)}+%
\overline{V_{i}})\cap X_{i}\subset U\right\} $

\ \ \ 

$U^{^{\prime }}$ is an open set, because $W_{i}$ is open, $\left\{ x\in W_{i}%
\text{ }\mid \text{cl}F_{i}^{V_{i}}(x)\subset U\right\} $ open since cl$%
F_{i}^{V_{i}}(x)$ is an upper semicontinuous map on $W_{i}$ and the set $%
\{x\in X\mid ($cl$(B_{i_{x}})_{b_{i}(x)}+$cl$V_{i})\cap X_{i}\subset U\}$ is
open since $($cl$(B_{i_{x}})_{b_{i}(x)}+$cl$V_{i})\cap X_{i}$ is u.s.c. Then 
$T_{i}^{V_{i}}$ is upper semicontinuous on $X$ and has closed convex values.

Define $T^{V}:X\rightarrow 2^{X}$ by $T^{V}(x):=\underset{i\in I}{\prod }%
T_{i}^{V_{i}}(x)$ for each $x\in X$.

$T^{V}$ is an upper semicontinuous correspondence and has also non-empty
convex closed values.

Since $X$ is a compact convex set, by Fan's fixed-point theorem [5], there
exists $x_{V}^{\ast }\in X$ such that $x_{V}^{\ast }\in T^{V}(x_{V}^{\ast })$%
, i.e., for each $i\in I$, $(x_{V}^{\ast })_{i}\in T_{i}^{V_{i}}(x_{V}^{\ast
})$. If $x_{V}^{\ast }\in W_{i},$ $(x_{V}^{\ast })_{i}\in $cl$%
F_{i}^{V_{i}}(x_{V}^{\ast })$, which is a contradiction.

Hence, $(x_{V}^{\ast })_{i}\in $cl$((B_{i_{x_{V}^{\ast
}}})_{b_{i}(x_{V}^{\ast })}+V_{i})\cap X_{i}$ and $(A_{i_{x_{V}^{\ast
}}})_{a_{i}(x_{V}^{\ast })}\cap (P_{i_{x_{V}^{\ast }}})_{p_{i}(x_{V}^{\ast
})}=\emptyset ,$ i.e. $x_{V}^{\ast }\in Q_{V}$ where

$Q_{V}=\cap _{i\in I}\{x\in X:$ $x_{i}\in $cl$((B_{i_{x}})_{b_{i}(x)}+V_{i})%
\cap X_{i}$ and $(A_{i_{x}})_{a_{i}(x)}\cap (P_{i_{x}})_{p_{i}(x)}=\emptyset
\}.$

Since $W_{i}$ is open, $Q_{V}$ is the intersection of non-empty closed sets,
therefore it is non-empty, closed in $X$.

We prove that the family $\{Q_{V}:V\in \underset{i\in I}{\tprod }\text{\ss }%
_{i}\}$ has the finite intersection property.

Let $\{V^{(1)},V^{(2)},...,V^{(n)}\}$ be any finite set of $\underset{i\in I}%
{\tprod \text{\ss }_{i}}$ and let $V^{(k)}=(V_{i}^{(k)})_{i\in I}$, $%
k=1,...n.$ For each $i\in I$, let $V_{i}=\underset{k=1}{\overset{n}{\cap }}%
V_{i}^{(k)}$, then $V_{i}\in \text{\ss }_{i};$ thus $V=(V_{i})_{i\in I}\in 
\underset{i\in I}{\tprod }\text{\ss }_{i}.$ Clearly $Q_{V}\subset \underset{%
k=1}{\overset{n}{\cap }}Q_{V^{(k)}}$ so that $\underset{k=1}{\overset{n}{%
\cap }}Q_{V^{(k)}}\neq \emptyset .$

Since $X$ is compact and the family $\{Q_{V}:V\in \underset{i\in I}{\tprod }%
\text{\ss }_{i}\}$ has the finite intersection property, we have that $\cap
\{Q_{V}:V\in \underset{i\in I}{\tprod }\text{\ss }_{i}\}\neq \emptyset .$
Take any $x^{\ast }\in \cap \{Q_{V}:V\in \underset{i\in I}{\tprod }\text{\ss 
}_{i}\},$ then for each $i\in I$ and each $V_{i}\in \text{\ss }_{i},$ $%
x_{i}^{\ast }\in $cl$((B_{i_{x^{\ast }}})_{b_{i}(x^{\ast })}+V_{i})\cap
X_{i} $ and $(A_{i_{x^{\ast }}})_{a_{i}(x^{\ast })}\cap (P_{i_{x^{\ast
}}})_{p_{i}(x^{\ast })}=\emptyset ;$ but then $x_{i}^{\ast }\in \overline{%
B_{i}^{\prime }}(x^{\ast })$ from Lemma 2\emph{\ }and $(A_{i_{x^{\ast
}}})_{a_{i}(x^{\ast })}\cap (P_{i_{x^{\ast }}})_{p_{i}(x^{\ast })}$for each $%
i\in I$ so that $x^{\ast }$ is an equilibrium point of $\Gamma $ in $X$. $%
\square $

\end{document}